\documentclass{article}
\usepackage{amssymb,amsmath,amsthm}
\newtheorem{theorem}{Theorem}

\newtheorem{lemma}{Lemma}
\newtheorem{proposition}{Proposition}

\newtheorem{remark}{Remark}
\date{}
\numberwithin{equation}{section} \numberwithin{theorem}{section}
\numberwithin{lemma}{section} \numberwithin{corollary}{section}
\numberwithin{remark}{section} \numberwithin{proposition}{section}
\numberwithin{definition}{section}
\begin{document}
% let's define a new macro
\newcommand{\n}{\noindent}
\newcommand{\vs}{\vskip}

\title{ Continuity of the Free Boundary in a Problem involving the $A$-Laplacian}

\author{S. Challal$^1$ and A. Lyaghfouri$^{2}$\\
\\
$^1$ Glendon College, York University\\
Department of Mathematics\\
Toronto, Ontario, Canada\\
\\
$^2$ American University of Ras Al Khaimah\\
Department of Mathematics and Natural Sciences\\
Ras Al Khaimah, UAE\\
}\maketitle

\vs 0.5cm
\begin{abstract}
In this paper we investigate a two dimensional free boundary problem involving the
A-Laplacian. We show that the free boundary is represented
locally by graphs of a family of continuous functions.
\end{abstract}

\vs 0.5cm

\n Key words : $A$-Laplacian, Free boundary, Continuity

\vs 0.5cm \n AMS 2000 Mathematics Subject Classification: 35R35; 35J60

 \vs 0.5cm
\section*{Introduction}\label{S:intro}

\vs 0.5cm \n The authors have considered the following problem in \cite{[CL9]}

\begin{equation*}(P)
\begin{cases}
& \text{ Find }\,\,(u, \chi) \in  W^{1,A}(\Omega)\times L^\infty
(\Omega)
\text{ such that}:\\
& (i)\quad  0\leqslant u\leqslant M, \quad 0\leqslant
\chi\leqslant 1 ,
\quad u(1-\chi) = 0 \,\,\text{ a.e.  in } \Omega\\
 & (ii)\quad \Delta_A u=-div(\chi H(x)) \quad\text{in }
 (W_0^{1,A}(\Omega))',
\end{cases}
\end{equation*}

\n where $\Omega$ is an open bounded domain of
$\mathbb{R}^2$, $x=(x_1,x_2)$, $M$ is a positive constant,
$\displaystyle{A(t)=\int_0^t a(s)ds}$, 
$\displaystyle{\Delta_A=div\Big({{a(|\nabla u|)}\over{|\nabla u|}}\nabla u\Big)}$ 
is the $A$-Laplacian, $a$ is  a  $C^1$ function from
$ [0,\infty)$ to $ [0,\infty)$ such that $a(0)=0$ and for some positive
constants $a_0, a_1$

\begin{equation}\label{0.1}
a_0\leqslant {{ta'(t)}\over{a(t)}}\leqslant a_1    \qquad \forall t>0.
\end{equation}

\n As a consequence of (0.1), we have the following
monotonicity inequality (see \cite{[CL3]})

\begin{equation}\label{0.2}
 \Big( { a(|\xi|)\over |\xi|}  \xi-  { a(|\zeta|)\over |\zeta|}
 \zeta\Big). (\xi-\zeta) >0 \qquad \forall \xi, \zeta\in
 \mathbb{R}^2\setminus\{ 0\},\quad \xi\neq\zeta.
\end{equation}

\n For examples of functions $a(t)$, we refer to \cite{[CLR]}. 

\vs 0,5cm\n $H=(H_1,H_2)$ is a vector function that
satisfies for some positive constants $\underline{h}, \overline{h}$
\begin{eqnarray}\label{0.3-5}
&& |H_1|\leqslant \overline{h},\qquad  0< \underline{h} \leqslant H_2\leqslant  \overline{h}
 \quad\text{ in }\Omega \\
&&  H\in C^{0,1}(\overline{\Omega})\\
&&  div(H)\geqslant 0 \quad\text{ a.e. in } \Omega\\
&&  div(H)\leqslant \overline{h} \quad\text{ a.e. in }.
\end{eqnarray}

\vs 0.2cm \n The Orlicz-Sobolev space $W^{1,A}(\Omega)$ is defined  by:

\begin{eqnarray*}
&&W^{1,A}(\Omega)=\left\{\,u\in L^A(\Omega)\,:\,|\nabla u|\in
L^A(\Omega)\,\right\},\quad\text{where}\\
&&L^A(\Omega)=\left\{\,u\in L^1(\Omega)\,:\,\int_\Omega
A(|u(x)|)dx<\infty\,\right\}
\end{eqnarray*}

\vs 0.2cm\n Throughout this paper, we shall denote by $B_r(x)$ a ball with center
$x$ and radius $r$.

\vs 0.2cm\n In \cite{[CL9]}, we showed that for any solution $(u, \chi)$, $u$ is Lipschitz 
continuous and that the free boundary is a union of graphs of a family of lower 
semi-continuous functions depending only on the vector function $H$. In this paper, we will 
show that these functions are actually continuous and that $\chi$ is the characteristic
function of the set $\{u>0\}$.

\vs 0.2cm\n Problem $(P)$ describes a variety of free boundary problems
including the lubrication problem \cite{[AC]}, and the dam
problem \cite{[ChiL1]}, \cite{[ChiL2]}, \cite{[CaL]},
\cite{[CL1]}, \cite{[CL2]}, \cite{[CL8]}, \cite{[L2]}, \cite{[L3]}, 
\cite{[L4]}, \cite{[LZ1]}, and \cite{[LZ2]}.
For a more general framework, we refer to \cite{[C]}, \cite{[Cha]}, \cite{[L]}, 
\cite{[CL3]}, \cite{[CL5]}, \cite{[CL6]}, \cite{[CL7]} and \cite{[CL9]}.
Regarding the problem with a Newman boundary condition, we refer 
the reader to \cite{[LS]} and  \cite{[S]}.

\vs 0,5cm

\section{The free boundary}\label{1}

\vs 0.5cm \n The free boundary is defined as the intersection between the two
sets $\{u=0\}$ and $\overline{\{u>0\}}$. 
When $H_1=0$ and $H_2$ is a constant function, it is easy to show as in \cite{[CL7]} that 
$\chi_{x_2}\leqslant 0$ in ${\cal D}'(\Omega)$ and that the free boundary
is the graph of a continuous function $x_2=\phi(x_1)$. When $H$ is not a constant
vector, we can show as in \cite{[CL5]} that
\begin{equation}\label{1}
div (\chi H) - \chi(\{u>0\}) div ( H) \leqslant 0 \qquad\text{
in }\quad {\cal D}^\prime(\Omega).
\end{equation}

\n Actually (1.1) can be obtained from $(P)ii)$ and by adapting the proof of Lemma 1.4.
As a consequence of (1.1), a weak monotonicity of the function  
$\chi$  (see \cite{[CL9]}) holds, which means that $\chi$ decreases along the orbits
of the following differential equation:

\[\displaystyle{(E(w,h))\,\,\left\{\begin{array}
{r@{\quad=\quad}l}
X^\prime (t,w,h)& H(X(t,w,h))\\
 X(0,w,h) & (w,h)\end{array}\right.}\]

\n where $h\in\pi_{x_2}(\Omega)$ and $w\in\pi_{x_1}(\Omega\cap \{x_2=h\})$,
$\pi_{x_1}$ and $\pi_{x_2}$ are respectively the orthogonal projections on
the $x_1$ and $x_2$ axes. We will denote by $X(.,w)$ the  maximal solution
of $E(w,h)$ defined on the interval $(\alpha_{-} (w), \alpha_{+}(w))$.
We know \cite{[CL8]} that the limits
$\displaystyle{\lim_{t\rightarrow \alpha_{-} (w)^+}X(t,w)}\in \partial\Omega\cap
\{x_2<h\}$ and
$\displaystyle{\lim_{t\rightarrow \alpha_{+} (w)^-}X(t,w)}\in\partial\Omega\cap \{x_2>h\}$ both exit,
which we shall denote respectively by
$X(\alpha_{-} (w),w)$ and $X(\alpha_{+} (w),w)$.
We shall also denote the orbit of $X(.,w)$ by $\gamma(w)$.

\vs0.5cm\n Now, we recall for the reader's convenience a few technical properties
and definitions established in \cite{[CL5]} and \cite{[CL9]}:

\vs0.3cm $\bullet$  $\alpha_{+}$ and $\alpha_{-}$ are uniformly bounded.

\vs0.2cm $\bullet$ For each $h\in\pi_{x_2}(\Omega)$, 
the mapping is one to one
 \begin{eqnarray*}
  &&T_h\, :  \,{D}_h \longrightarrow  T_h({D}_h)\\
   &&\qquad(t,w)\longmapsto  T_h(t,w)=(T_h^1,T_h^2)(t,w)=X(t,w).
 \end{eqnarray*}
where $\displaystyle{D_h=\{ (t,w)\,/\, w\in \pi_{x_1}(\Omega\cap\{x_2=h\}), \, t\in
(\alpha_{-}(w), \alpha_{+}(w))\}}$.

\vs0.2cm $\bullet$  $\displaystyle{ \Omega=\bigsqcup_{h\in\pi_{x_2}(\Omega)}T_h(D_h)}$.

\vs0.2cm $\bullet$ $ T_h$  and $T_h^{-1} $ are $C^{0,1}$.
\vs 0,2cm $\bullet$  The determinant of the Jacobian  matrix  of the
mapping $T_h$, denoted by $Y_h(t,w)$,  satisfies:

\vs 0.2cm\n i) \quad $\displaystyle{ Y_h(t,w)=
-H_2(w,h) \exp\big( \int_0^t  (div H)(X(s,w)) ds\big)}
\qquad\mbox{a.e. in } D_h.$

\vs 0.2cm \n ii) \quad  $\underline{h} \leqslant
-Y_h(t, w) \leqslant C {\overline{h}},\quad C>0,
\qquad\mbox{a.e. in } D_h.$

\vs 0.5cm \n The following monotonicity of
$\chi$ based on (1.1) (see \cite{[CL5]},  \cite{[CL9]}) is the key point in
parameterizing the free boundary:

\begin{equation}\label{2}
{\partial \over\partial t} \big(\chi o T_h \big)\leqslant 0\qquad\text{ in }
{\cal D}'(D_h).
\end{equation}

\vs 0.5cm\n Property (1.2) means that $\chi$ decreases along the orbits
of the differential equation $(E(w,h))$. The consequence
of this monotonicity on $u$ is materialized in the next theorem
established in \cite{[CL9]}.

\vs 0,3cm
\begin{theorem}\label{t2.1.} Let $(u,\chi)$ be a solution of
$(P)$ and $x_0=T_h(t_0,w_0)\in T_h(D_h)$.

\vs 0,2cm \n
 i)\quad If $\quad u(x_0)=uoT_h(t_0,w_0) > 0,\quad$ then there exists $\epsilon >0$
  such that
$$uoT_h(t,w) > 0 \qquad \forall
 (t,w)\in C_\epsilon=\{ (t,w)\in D_h\,/ \,\, \vert w-
w_0\vert <\epsilon , \,\, t <t_0 +\epsilon \}.$$
\vs 0,2cm
\n ii)\quad If $\quad u(x_0)=uoT_h(t_0,w_0)=0,\quad $
then $\quad uoT_h(t, w_0) =0 \qquad \forall t\geqslant
t_0.$\end{theorem}

\vs 0.3cm\n The proof of Theorem 1.1 is based on the following
strong maximum principle:
\vs 0,3cm

\begin{lemma}\label{l1.1}

\n Let $u\in W^{1,A}(U)\cap C^1(U)\cap C^0(\overline{U})$ such that $ u
\geqslant 0$ in $U$ and $\Delta_A u\leqslant 0$ in $U$.
Then $u\equiv 0$ in $U$ or $u>0$ in $U$.
\end{lemma}

\vs 0,3cm \n Thanks to Theorem 1.1, we can define for each
$h\in \pi_{x_2}(\Omega)$, the following function $\phi_h$ on
$\pi_{x_1}(\Omega\cap \{x_2=h\})$ (see \cite{[CL9]}):

\begin{equation*}\phi_h(w) =\begin{cases}
&\sup\left\{t\,~:~ \, (t,w)\in D_h, \quad u o T_h(t,w) > 0\right\}
 \\
&\hskip 2cm \text{ if this set is not empty }\\
 & \alpha_-(w)\hskip 2cm
\text{ otherwise.}\end{cases}
\end{equation*}

\n Then one can easily establish the following proposition  as in \cite{[CL1]}:

\vs 0,3cm
\begin{proposition}\label{p1.1.} For each $h\in \pi_{x_2}(\Omega)$, the
function $\phi_h$ is lower semi-continuous at each $w \in\pi_{x_1}(\Omega\cap\{x_2=h\})$
such that $T_h(\phi_h(w), w) \in\Omega$. Moreover
\begin{equation}\label{e1.3.}
\{u oT_h(t,w)
>0\}\cap D_h=\{t<\phi_h(w)\}.
\end{equation}
\end{proposition}

\vs 0,3cm
\begin{lemma}\label{l1.2}
Let $h\in \pi_{x_2}(\Omega)$. For each $k\in\pi_{x_2}(\Omega)$ and $w\in\pi_{x_1}(\Omega\cap \{x_2=h\})$,
let $t_k(w)$ be the unique value of $t$  at which the orbit
$\gamma(w)$ intersects with the line $\{x_2=k\}$ if it exists. Then the function
$S(k,w)=t_k(w)$ is Lipschitz continuous in its domain. More precisely, we have
for some positive constant $C$
\begin{eqnarray*}
|S(k,w)-S(k_0,w_0)|\leqslant C(|k-k_0|+|w-w_0|)\quad
\forall (k,w), (k_0,w_0)\in domain(S).
 \end{eqnarray*}
\end{lemma}

\vs0.3cm\n \emph{Proof.} Let $(k,w), (k_0,w_0)\in domain(S)$. First we have from the differential 
equation $(E(w,h))$
\begin{eqnarray*}
k=h+\int_0^{t_k(w)} H_2(X(t,w))ds\quad\text{ and }\quad
k_0=h+\int_0^{t_{k_0}(w_0)} H_2(X(t,w_0))ds.
\end{eqnarray*}

\n If we subtract the two equalities from each other, we get
\begin{eqnarray}\label{e1.4}
&& k-k_0=\int_0^{t_k(w)} H_2(X(t,w))ds-\int_0^{t_{k_0}(w_0)} H_2(X(t,w_0))ds.
\end{eqnarray}

\n Next, if we assume that $t_k(w)>t_{k_0}(w_0)$,
then we get by (0.3)
\begin{eqnarray}\label{e1.5}
\underline{h}(t_k(w)-t_{k_0}(w_0))&\leqslant& \int_{t_{k_0}(w_0)}^{t_k( w)}H_2(X(s,w))ds).
\end{eqnarray}

\n Now, observe that
\begin{eqnarray}\label{e1.6}
\int_{t_{k_0}(w_0)}^{t_k( w)}H_2(X(s,w))ds)&=&\int_0^{t_k(w)} H_2(X(t,w))ds-\int_0^{t_{k_0}(w_0)} H_2(X(t,w))ds\nonumber\\
&=&\int_0^{t_k(w)} H_2(X(t,w))ds-\int_0^{t_{k_0}(w_0)} H_2(X(t,w_0))ds\nonumber\\
&&+\int_0^{t_{k_0}(w_0)} (H_2(X(t,w_0))ds-H_2(X(t,w)))ds.
\end{eqnarray}

\n Using (1.4), (1.6) and the fact that $H_2\circ X$ is Lipschitz continuous in $\overline{D_h}$, and
since $t_{k_0}(w_0)$ is bounded independently of $k$ and $w$, we obtain from (1.5), for some positive 
constant $C_0$

\begin{eqnarray*}
\underline{h}(t_k(w)-t_{k_0}(w_0))&\leqslant& k-k_0+ C_0|w-w_0|
\end{eqnarray*}

\n which leads for $C=\max(1,C_0)/\underline{h}$, to
\begin{eqnarray}\label{e1.7}
t_k(w)-t_{k_0}(w_0)&\leqslant& C(|k-k_0|+|w-w_0|).
\end{eqnarray}

\n If $t_k(w)<t_{k_0}(w_0)$, we get in a similar fashion
\begin{eqnarray}\label{e1.8}
t_{k_0}(w_0)-t_k(w)&\leqslant& C(|k-k_0|+|w-w_0|).
\end{eqnarray}

\n Combining (1.7) and (1.8), the lemma follows.
\qed
\vs 0.3cm\n Our main goal is to prove that for each $h\in \pi_{x_2}(\Omega)$, the
function $\phi_h$ is actually continuous. Due to the
local character of this result, we will confine ourselves to the following
situation:

\vs 0.2cm \n We assume that $u=0$ on an open and connected subset $\Gamma$ of
$\partial\Omega$ and consider a free boundary point in
a neighborhood of the form $U= T_h \big((w_*, w^*)\times
(h,\infty)\big) \cap\Omega \subset T_h
(D_h) $ such that $T_h\left( \left\{(w, \alpha_+(w)),\,\,
w\in (w_*, w^*)\right\}\right) \subset\subset\Gamma$.
So we are led to study the following problem:

$$ (P)\left\{
\begin{array}{ll}
\text{Find } (u,\chi)\in  W^{1,A}(U)\times L^\infty(U)
\text{ such that }: \quad u=0\hbox{ on } \Gamma& \\
\quad  0\leqslant u\leqslant M, \quad 0\leqslant
\chi\leqslant 1 ,
\quad u(1-\chi) = 0 \,\,\text{ a.e.  in } U & \hbox{} \\
\displaystyle{\int_{U}}\Big( {a(|\nabla u|)\over |\nabla u|} \nabla u +\chi H(x)\Big)
.\nabla\zeta dx \leqslant 0& \hbox{} \\
\quad\forall\zeta\in W^{1,A}(U),\quad \zeta\geqslant 0 \quad \hbox{ on } \Gamma, \quad
\zeta  =0 \quad\hbox{ on } \partial U\setminus \Gamma.& \hbox{}
  \end{array}
\right.
$$

\vs 0.5cm\n  We observe that the free boundary $(\partial\{u>0\})\cap U$ is
the graph of the lower semi-continuous function $\phi_h$ in $(w_*, w^*)$.
Our objective is to prove the continuity of the function $\phi_h$. To do
that, it is enough to show that it is upper semi-continuous. To this end, we need to
generalize few lemmas previously established for a linear operator in \cite{[CL5]}.
In the sequel and without notice, we shall denote by $(u,\chi)$ a solution of 
the problem $(P)$.

\vs 0,3cm
\begin{lemma}\label{l1.3}
Let $w_1, w_2 \in(w_*, w^*)$  with $w_1< w_2$,
and let  $k \in \pi_{x_2}(U)$ such that
$\{x_2=k\} \cap \gamma(w_i) \neq\emptyset$,  $ i=1, 2$.
Let
$$Z_k=  T_h\left(\left\{
(t,w)\in D_h, \,\, w\in (w_1, w_2),\,\,
t>t_k(w)\right\}\right)= T_h(\{w_1<
w<w_2\})\cap \left\{x_2> k\right\}.$$
If   $uoT_h(t_k(w_i),w_i) =0 $ for $i=1,2$, then we have

\begin{eqnarray*}
&&\int_{Z_k}\Big( {a(|\nabla u|)\over |\nabla u|} \nabla u +\chi H(x)\Big)
.\nabla\zeta dx \leqslant 0\\
&&\quad\forall\zeta\in W^{1,A}(Z_k)\cap C^0(\overline{Z}_k),\quad \zeta\geqslant 0
\quad\hbox{on}\quad\overline{Z}_k\setminus \{x_2=k\},\quad
\zeta=0 \quad\hbox{on}\quad\overline{Z}_k\cap \{x_2=k\}.
 \end{eqnarray*}

\end{lemma}

\vs0,3cm\n The proof of Lemma 1.3 is inspired from the one of a similar lemma 
in \cite{[C]} obtained for the case $H(x)=(h(x),0)$. Our proof is based 
on the next lemma.

\begin{lemma}\label{l1.4} Under the assumptions of Lemma 1.3, we have
\begin{eqnarray*}
&&\int_{Z_k} \Big( {a(|\nabla u|)\over |\nabla u|} \nabla u.\nabla\zeta dx -\int_{Z_k}\chi_{\{u>0\}}div(H).\zeta dx \leqslant 0\\
&&\quad\forall\zeta\in W^{1,A}(Z_k)\cap C^0(\overline{Z}_k),\quad \zeta\geqslant 0
\quad\hbox{on}\quad\overline{Z}_k\setminus \{x_2=k\},\quad
\zeta=0 \quad\hbox{on}\quad\overline{Z}_k\cap \{x_2=k\}.
\end{eqnarray*}
\end{lemma}

\n \emph{Proof.} Let $\zeta$ be as in the lemma,
$\epsilon>0$, and $F_\epsilon(u)=\min\big({u^+\over\epsilon},1\big)$.
Using $\chi(Z_k)F_\epsilon(u)\zeta$ as a test function for $(P)$, we get
\begin{eqnarray*}
&&\int_{Z_k}F_\epsilon(u) {a(|\nabla u|)\over |\nabla u|} \nabla u.\nabla\zeta dx
+\int_{Z_k} H(x).\nabla(F_\epsilon(u)\zeta) dx\\
&&~~\leqslant-\int_{Z_k}
F_\epsilon'(u)\zeta  |\nabla u|a(|\nabla u|) dx\leqslant0\\
\end{eqnarray*}

\n Integrating by parts, we obtain
\begin{eqnarray}\label{e1.9}
\int_{Z_k}F_\epsilon(u) {a(|\nabla u|)\over |\nabla u|} \nabla u\nabla\zeta dx
-\int_{Z_k}div(H).F_\epsilon(u)\zeta dx \leqslant 0.
\end{eqnarray}

\n The lemma follows by letting $\epsilon$ go to $0$ in (1.9).
\qed

\vs0,3cm\n\emph{Proof.} For $\epsilon>0$ small enough, let
$\alpha_\epsilon(w)
=\displaystyle\min\left(1,\frac{(w-w_1)^+}{\epsilon},\frac{(w_2-w)^+}{\epsilon}\right)$,
and observe that
\begin{eqnarray}\label{e1.10}
&&\int_{Z_k}\Big( {a(|\nabla u|)\over |\nabla u|} \nabla u + \chi
H(x)\Big) .\nabla \zeta dx=\int_{Z_k} \Big( {a(|\nabla u|)\over |\nabla u|} \nabla u + \chi
H(x)\big) .\nabla[(\alpha_\epsilon\circ T_h^{-1})\zeta] dx\nonumber\\
&&+\int_{Z_k} \Big( {a(|\nabla u|)\over |\nabla u|} \nabla u + \chi
H(x)\Big) .\nabla[(1-\alpha_\epsilon\circ T_h^{-1})\zeta] dx
\end{eqnarray}

\n Since $\chi(Z_k)(\alpha_\epsilon\circ T_h^{-1})\zeta$ is a
test function for (P), we have:
\begin{equation}\label{e1.11}
\displaystyle{\int_{Z_k} }\Big( {a(|\nabla u|)\over |\nabla u|} \nabla u + \chi
H(x)\Big) .\nabla[(\alpha_\epsilon\circ T_h^{-1})\zeta] dx
\,\leqslant\,0
\end{equation}

\n Now, if we apply Lemma 1.4 to the function $(1-\alpha_\epsilon\circ T_h^{-1})\zeta$,
we get
\begin{eqnarray}\label{e1.12}
&&\int_{Z_k} {a(|\nabla u|)\over |\nabla u|} \nabla u.\nabla[(1-\alpha_\epsilon\circ T_h^{-1})\zeta] dx \nonumber\\
&&\leqslant \int_{Z_k}\chi_{\{u>0\}}div(H).(1-\alpha_\epsilon\circ T_h^{-1})\zeta dx 
\end{eqnarray}

\n Taking into account (1.11)-(1.12), we obtain from (1.10)
\begin{eqnarray}\label{e1.13}
&&\int_{Z_k}\Big( {a(|\nabla u|)\over |\nabla u|} \nabla u  + \chi
H(x)\big) .\nabla \zeta dx\nonumber\\
&&\leqslant \int_{Z_k}\chi_{\{u>0\}}div(H).(1-\alpha_\epsilon\circ T_h^{-1})\zeta dx
+\int_{Z_k} \chi H(x).\nabla[(1-\alpha_\epsilon\circ T_h^{-1})\zeta] dx
\end{eqnarray}

\n Using the change of variables $x=T_h(t,w)$ and arguing as in the proof of Theorem 2.1
in \cite{[CL5]}, we get
\begin{eqnarray}\label{e1.14}
\int_{Z_k} \chi H(x).\nabla[(1-\alpha_\epsilon\circ T_h^{-1})\zeta] dx
&=&\int_{T_h^{-1}(Z_k)} Y_h \chi\circ T_h
\frac{\partial}{\partial t}[1-\alpha_\epsilon\zeta\circ T_h]dtdw\nonumber\\
&=&\int_{T_h^{-1}(Z_k)} (1-\alpha_\epsilon)Y_h \chi\circ T_h
\frac{\partial}{\partial t}[\zeta\circ T_h]dtdw
\end{eqnarray}

\n Hence we derive from (1.13) and (1.14)
\begin{eqnarray}\label{e1.15}
&&\int_{Z_k}\Big( {a(|\nabla u|)\over |\nabla u|} \nabla u + \chi
H(x)\Big) .\nabla \zeta dx  \nonumber\\
&&\leqslant\int_{Z_k}\chi_{\{u>0\}}div(H).(1-\alpha_\epsilon\circ T_h^{-1})\zeta dx\nonumber\\
&&+\int_{T_h^{-1}(Z_k)} (1-\alpha_\epsilon)Y_h \chi\circ T_h
\frac{\partial}{\partial t}[\zeta\circ T_h]dtdw
\end{eqnarray}

\n The lemma follows by letting $\epsilon$ go to $0$ in (1.15).
\qed

\vs 0,3cm
\begin{lemma}\label{l1.2}
Let $x_0=T_h(t_0,w_0)$ be a point in $U$.

\n If $u oT_h= 0 $ in $B_r(t_0,w_0)$, then
$$u oT_h= 0 \quad\hbox{in }  C_r
\quad\hbox{ and }\quad \chi oT_h= 0 \quad\hbox{a.e. in } C_r$$
where $C_r=\{ (t,w)\in D_h, \quad \vert w -w_0\vert
< r, \quad t >t_0\} \cup B_r (t_0,w_0).$
\end{lemma}

\vs 0,3cm\n \emph{Proof.} By Theorem 1.1 $ii)$, we have  $uoT_h=0$ in
$C_r$. Applying  Lemma 1.2 with  domains $Z_k=T_h(\{
w_1<w<w_2\})\cap \{x_2>k\} \subset
T_h(C_r)$, ($k\in\pi_{x_2}(U)$) and taking
$\zeta=x_2-k$, we obtain
$\displaystyle{\int_{Z_k}\chi H_2 dx \leqslant 0}$.
Then we deduce from (0.3) that  $\chi=0$ a.e. in
$Z_k$. This holds for all domains
$Z_k$ in $T_h(C_r)$. Hence  $\chi =0$ a.e. in
$T_h(C_r)$.\qed

\vs 0,3cm

\begin{lemma}\label{l1.3}
Let $x_0=T_h(t_0,w_0)\in U$ such that $B_r=B_r(t_0,w_0)\subset D_h$.
Then we cannot have the following three situations
\begin{eqnarray*}
&(i) &\begin{cases}
uoT_h(t,w_0)=0\qquad\forall t\in(t_0-r,t_0+r)\\
u o T_h(t,w)>0\qquad\forall t\in(t_0-r,t_0+r),\quad \forall w\neq w_0,
\end{cases}\\
&(ii)& \begin{cases} u o T_h(t,w)=0\qquad\forall(t,w)\in
B_r\cap \{w\leqslant w_0\}\\
u o T_h(t,w)>0\qquad\forall(t,w)\in B_r\cap \{w>w_0\},\end{cases}\\
&(iii)& \begin{cases}
 u o T_h(t,w)=0\qquad\forall(t,w)\in B_r\cap\{w\geqslant w_\}] \\
u o T_h(t,w)>0\qquad\forall(t,w)\in B_r\cap \{w<w_0\}.
\end{cases}\end{eqnarray*}
\end{lemma}

\vs0.3cm\n \emph{Proof.}  Assume that $ii)$ holds. The proofs of $i)$ and $iii)$ are based
on similar arguments. Let  $\zeta\in\mathcal{D}(T_h(B_r))$,
$\zeta\geqslant 0$. Using the fact that, by Lemma 1.3, $\chi oT_h
=0$ a.e. in $B_r \cap \{w\leqslant w_0\}$, we obtain after using
the change of variable $T_h$

$$\int_{T_h (B_r)} {a(|\nabla u|)\over |\nabla u|} \nabla u . \nabla\zeta dx = \int_{B_r\cap [w>w_0]}
{\partial\over \partial t} (-Y_h(t,w)) \zeta  oT_h dt
dw\geqslant 0.$$

\n This means that $\triangle_A u \leqslant 0$ in $\mathcal{D}^{\prime}(T_h(B_r))$.
By Lemma 1.1, either $u>0$ or $u=0$ in $T_h(B_r)$, which contradicts
the assumption.\qed

\vs 0.5cm
\section{Continuity of the free boundary}\label{S3}

\vs 0.5cm \n As mentioned in section 1, to prove the continuity of the function $\phi_h$, 
it is enough to show that it is upper semi-continuous. The main idea to do that is to compare $u$
with a suitable barrier function near a free boundary point. In the following  step, we construct
such a function. For this purpose, let $\epsilon>0$, $w_1, w_2 \in(w_*, w^*)$
such that $w_1< w_2$, $k \in \pi_{x_2}(U)$,
and assume that $\epsilon$ is small enough to guarantee that 
$Z_k^{k+\epsilon}(w_1,w_2)=T_h(\{w_1<w<w_2\})\cap\{k<x_2<k+\epsilon\}\subset\subset U$
and $\epsilon<\underline{h}/{2\overline{h}}$.

\n Then the function $\overline{v}_\epsilon$ defined by

$$\overline{v}_\epsilon(x_1,x_2)=\vartheta_\epsilon(k+\epsilon -x_2)  \quad\hbox{ with }\quad
 \vartheta_\epsilon (t)=\int_{0}^{t} a^{-1}(2 \overline{{h}}\epsilon
- \overline{{h}}s) ds \quad\text{ for } t\in[0,\epsilon]$$

\n satisfies

\begin{equation}\label{2.1}
\Delta_A \overline{v}_\epsilon= - \overline{h} \quad \text{in  }Z_k^{k+\epsilon}(w_1,w_2).
\end{equation}

\n Now let $v_\epsilon$ be the unique solution in $W^{1,A}(Z_k^{k+\epsilon}(w_1,w_2))$ of

\begin{equation}\label{2.2}
  \begin{cases}
   \quad \Delta_A v_\epsilon= - div(H) & \quad \text{in  }Z_k^{k+\epsilon}(w_1,w_2) \\
   \quad v_\epsilon=\overline{v}_\epsilon &\quad  \text{on  }\partial Z_k^{k+\epsilon}(w_1,w_2).
  \end{cases}
\end{equation}

\n Then we have:
\vs 0,3cm

\begin{lemma}\label{l2.1}

\begin{equation}\label{2.3}
0\leqslant v_\epsilon \leqslant \overline{v}_\epsilon
 \quad \hbox{ in }\quad Z_k^{k+\epsilon}(w_1,w_2).
\end{equation}

\end{lemma}

\vs0,3cm\n\emph{Proof.} $i)$ Note that $v_\epsilon^-\in W^{1,A}(Z_k^{k+\epsilon}(w_1,w_2))$ and
$v_\epsilon^-=0$ on $\partial Z_k^{k+\epsilon}(w_1,w_2)$. Therefore we obtain from (2.2) and (0.5)

\begin{eqnarray}\label{2.4}
&& \int_{Z_k^{k+\epsilon}(w_1,w_2)}{{a(|\nabla v_\epsilon|)}\over{|\nabla
v_\epsilon|}}\nabla v_\epsilon  .\nabla v_\epsilon^- dx\,=\,
\int_{Z_k^{k+\epsilon}(w_1,w_2)} div(H) .v_\epsilon^- dx \nonumber\\
&&\int_{Z_k^{k+\epsilon}(w_1,w_2)}{{a(|\nabla v^-|)}\over{|\nabla
v_\epsilon^-|}}\nabla v_\epsilon^-  .\nabla v_\epsilon^- dx\,=\,
\int_{Z_k^{k+\epsilon}(w_1,w_2)} -div(H) .v_\epsilon^- dx\nonumber\\
&&\int_{Z_k^{k+\epsilon}(w_1,w_2)}|\nabla v_\epsilon^-|a(|\nabla v_\epsilon^-|)dx \,=\,
\int_{Z_k^{k+\epsilon}(w_1,w_2)} -div(H) .v_\epsilon^- dx \leqslant 0.
\end{eqnarray}

\n Taking into account (2.4) and the fact that $ta(t)$ is an
increasing function, we deduce that $\nabla v_\epsilon^-=0$ a.e. in $Z_k^{k+\epsilon}(w_1,w_2)$.
Since $v_\epsilon^-=0$ on $\partial Z_k^{k+\epsilon}(w_1,w_2)$, we must have $v_\epsilon^-=0$ in $Z_k^{k+\epsilon}(w_1,w_2)$.
Hence $v_\epsilon\geqslant 0$ in $Z_k^{k+\epsilon}(w_1,w_2)$.

\vs0,2cm\n $ii)$ Similarly, we observe that $(v_\epsilon-\overline{v}_\epsilon)^+\in W^{1,A}(Z_k^{k+\epsilon}(w_1,w_2))$ and
$(v_\epsilon-\overline{v}_\epsilon)^+=0$ on $\partial Z_k^{k+\epsilon}(w_1,w_2)$. Therefore we obtain
from (2.1) and (2.2)

\begin{eqnarray}\label{2.5-6}
&& \int_{Z_k^{k+\epsilon}(w_1,w_2)}{{a(|\nabla v_\epsilon|)}\over{|\nabla
v_\epsilon|}}\nabla v_\epsilon  .\nabla (v_\epsilon-\overline{v}_\epsilon)^+ dx\,=\,
\int_{Z_k^{k+\epsilon}(w_1,w_2)} div(H) .(v_\epsilon-\overline{v}_\epsilon)^+ dx \\
&&\int_{Z_k^{k+\epsilon}(w_1,w_2)}{{a(|\nabla \overline{v}_\epsilon|)}\over{|\nabla
\overline{v}_\epsilon|}}\nabla \overline{v}_\epsilon  .\nabla (v_\epsilon-\overline{v}_\epsilon)^+ dx\,=\,
\int_{Z_k^{k+\epsilon}(w_1,w_2)} \overline{h} .(v_\epsilon-\overline{v}_\epsilon)^+ dx.
\end{eqnarray}

\n Subtracting (2.6) from (2.5), and using (0.6), we get
\begin{eqnarray}\label{2.4}
&& \int_{Z_k^{k+\epsilon}(w_1,w_2)}\left({{a(|\nabla v_\epsilon|)}\over{|\nabla
v_\epsilon|}}\nabla v_\epsilon-{{a(|\nabla \overline{v}_\epsilon|)}\over{|\nabla
\overline{v}_\epsilon|}}\nabla \overline{v}_\epsilon\right)  .\nabla (v_\epsilon-\overline{v}_\epsilon)^+dx\nonumber\\
&&\quad=\int_{Z_k^{k+\epsilon}(w_1,w_2)} (div(H)-\overline{h}) .(v_\epsilon-\overline{v}_\epsilon)^+ dx\leqslant 0.
\end{eqnarray}

\n Taking into account (2.7) and (0.2), we obtain $\nabla (v_\epsilon-v_\epsilon)^+=0$ a.e. in $Z_k^{k+\epsilon}(w_1,w_2)$.
Since $(v_\epsilon-\overline{v}_\epsilon)^+=0$ on $\partial Z_k^{k+\epsilon}(w_1,w_2)$, we get $(v_\epsilon-\overline{v}_\epsilon)^+=0$ in $Z_k^{k+\epsilon}(w_1,w_2)$.
Hence $v_\epsilon\leqslant \overline{v}_\epsilon$ in $Z_k^{k+\epsilon}(w_1,w_2)$.

\qed

\vs 0,3cm
\begin{lemma}\label{l2.2} 
After extending $v_\epsilon$ by $0$  to
$Z_{k+\epsilon}$, we obtain
$$ \int_{Z_k}\Big( {a(|\nabla v_\epsilon|)\over |\nabla v_\epsilon|} \nabla v_\epsilon + \chi([v>0]) H(x)\Big)\nabla\zeta dx \geqslant 0
\quad \forall \zeta\in W^{1,A}(Z_k),\, \zeta\geqslant 0, \,\zeta=0 \hbox{ on } \partial Z_k\cap U.
 $$
\end{lemma}

\vs 0,3cm\n\emph{Proof.} First we have $\Delta_A v_\epsilon= - div H\leqslant 0$ in $Z_k^{k+\epsilon}(w_1,w_2)$,
and by (2.3), $v_\epsilon\geqslant 0$ in $Z_k^{k+\epsilon}(w_1,w_2)$. By Lemma 1.1, we obtain $v_\epsilon>0$ in $Z_k^{k+\epsilon}(w_1,w_2)$.

\n Note that $v_\epsilon=0$ on $L=\partial Z_k^{k+\epsilon}(w_1,w_2)\cap\{x_2=k+\epsilon\}$ and
$v_\epsilon\in C^{1, \alpha}_{loc}(Z_k^{k+\epsilon}(w_1,w_2)\cup L)$  for some $\alpha\in (0,1)$ (see \cite{[L]}).

\vs0,2cm\n  Next we claim that

\begin{equation}\label{2.8}
|\nabla v_\epsilon(x)|\leqslant
a^{-1}(2\overline{h}\epsilon)\qquad\forall x\in L.
\end{equation}

\n From Lemma 2.1, we have $v_\epsilon\leqslant \overline{v}_\epsilon$ in $Z_k^{k+\epsilon}(w_1,w_2)$.
In particular we have
\begin{eqnarray*}
\forall (x_1,x_2)\in Z_k^{k+\epsilon}(w_1,w_2)\qquad v_\epsilon(x_1,x_2)-v(x_1,k+\epsilon)
&\leqslant&\overline{v}_\epsilon(x_1,x_2)-\overline{v}_\epsilon(x_1,k+\epsilon)\\
&\leqslant&\max_{\overline{Z}_k(w_1,w_2)}|\nabla \overline{v}_\epsilon||x_2-k-\epsilon|.
\end{eqnarray*}

\n We obtain  $(2.8)$ since
$$|\nabla \overline{v}_\epsilon|=\vartheta_\epsilon'(k+\epsilon-x_2)\leqslant
\vartheta_\epsilon'(0)=a^{-1}(2\overline{h}\epsilon).$$

\n Now if $\nu$ is the outward unit normal vector to $L$, then we have by (2.8),
since $\epsilon\in (0, \underline{h}/{2\overline{h}})$
\begin{eqnarray}\label{2.9}
{{a(|\nabla v_\epsilon|)}\over{|\nabla v_\epsilon|}} \nabla v_\epsilon.\nu +
H(x).\nu &=&{a(|\nabla v_\epsilon|)\over |\nabla v_\epsilon|} \nabla v_\epsilon.e_2 + H_2(x)
\geqslant - a(|\nabla v_\epsilon|)+ \underline{h} \nonumber\\
&\geqslant& -2\overline{h}\epsilon + \underline{h}\geqslant 0\quad\text{on}~~L
\end{eqnarray}

\n Finally, for $\zeta\in W^{1,A}(Z_k)$, $\zeta\geqslant 0$,
$\zeta=0$ on $\partial Z_k\cap U$, we obtain from (2.2) and (2.9)
$$ \int_{Z_k}\Big( {a(|\nabla v_\epsilon|)\over |\nabla v_\epsilon|} \nabla v_\epsilon + \chi(\{v_\epsilon>0\}) H(x)\Big)\nabla\zeta dx=
\int_{L}\Big( {a(|\nabla v_\epsilon|)\over |\nabla v_\epsilon|}  \nabla v_\epsilon .\nu+  H(x).\nu\Big)\zeta dx \geqslant
0.$$
\qed

\vs0,3cm

\begin{lemma}\label{l2.2}
Assume that 
\begin{eqnarray*}
  && uoT_h (t_k(w_1), w_1)
  =uoT_h (t_k(w_2), w_2)= 0\\
  && uoT_h (t_k(w), w)\leqslant \vartheta_\epsilon(\epsilon)=
  v_\epsilon (t_k(w), w)\qquad \forall w\in
  (w_1, w_2).
\end{eqnarray*}
Then we have
$$\lim_{\delta\rightarrow 0}{1\over \delta}\int_{ Z_k^{k+\epsilon}(w_1,w_2)\cap
  \{0< u-v_\epsilon< \delta\}} \Big( {a(|\nabla u|)\over |\nabla u|}  \nabla u-
  {a(|\nabla v_\epsilon|)\over |\nabla v_\epsilon|}\nabla v_\epsilon\Big).  \nabla (u-v_\epsilon)  dx  = 0
$$
\end{lemma}

\vs 0,3cm\n \emph{Proof.}  For $\delta,\eta> 0$, let $F_\delta (s)$ as in the proof of Lemma 1.4, 
$\, d_\eta (x_2)=F_\eta (x_2-\overline{k})$  and $\overline{k}= k+\epsilon$.
By applying Lemma 1.2 and Lemma 2.1 respectively for $\zeta = F_\delta(u-v_\epsilon) + d_\eta (1-F_\delta (u))$ and
$\zeta = F_\delta(u-v_\epsilon)$, we get
\begin{eqnarray*}
&&\int_{Z_k}\Big( {a(|\nabla u|)\over |\nabla u|}
\nabla u + \chi H(x)\Big).\nabla (F_\delta (u-v_\epsilon) ) dx \nonumber\\
&&\qquad \leqslant - \int_{Z_k}\Big( {a(|\nabla
u|)\over |\nabla u|}   \nabla u + \chi H(x) \Big).\nabla (d_\eta(1- F_\delta (u)) ) dx
\end{eqnarray*}

\begin{equation*}
-\int_{Z_k}\Big( {a(|\nabla v_\epsilon|)\over |\nabla v_\epsilon|}  \nabla v_\epsilon + \chi(\{v_\epsilon>0\}) H(x)
  \Big).\nabla (F_\delta (u-v_\epsilon) ) dx\leqslant 0.
\end{equation*}
Adding these inequalities, we get  since $d_\eta=0$ in $\{v_\epsilon>0\}$

\begin{eqnarray*} &&
\int_{Z_k\cap\{v_\epsilon>0\}} F_\delta ^{\prime}(u-v_\epsilon) \Big(
{a(|\nabla u|)\over |\nabla u|}  \nabla u - {a(|\nabla v_\epsilon|)\over
|\nabla v_\epsilon|} \nabla v_\epsilon\Big) .
\nabla (u-v_\epsilon) dx\\
&&\quad \leqslant  \displaystyle{-\int_{Z_k\cap
\{v_\epsilon=0\}} (1- d_\eta)\Big( {a(|\nabla u|)\over |\nabla u|}   \nabla u
+ \chi H(x)\Big).\nabla ( F_\delta (u) ) dx}\\
 &&\qquad -  \displaystyle{
 \int_{Z_k\cap\{v_\epsilon=0\}} (1- F_\delta (u))\Big( {a(|\nabla u|)\over |\nabla u|}   \nabla u + \chi H(x)
 \Big).\nabla d_\eta dx  = I_1^{\delta\eta} + I_2^{\delta\eta}.}
\end{eqnarray*}

\n Since
$$|I_1^{\delta\eta}|\leqslant \displaystyle{\int_{D_{k\cap\{\overline{k}<x_2<\overline{k}+\eta\}}}
\big(a(|  \nabla u|) + | H(x)| \big)|\nabla ( F_\delta (u) )|
dx},$$ we obtain $\displaystyle{ \lim_{\eta\rightarrow 0}
I_1^{\delta\eta}=0}$.

\n As for $I_2^{\delta\eta}$, we have

\begin{eqnarray*}  &&I_2^{\delta\eta}= -\int_{Z_k
\cap [u=v_\epsilon=0]} \chi H(x).\nabla d_\eta dx\\
&&\qquad- \int_{Z_k\cap [u>v_\epsilon=0]} (1- F_\delta
(u))\Big( {a(|\nabla u|)\over |\nabla u|}  \nabla u + H(x)
\Big).\nabla d_\eta dx  \leqslant I_3^{\delta\eta}
\end{eqnarray*}

\n since we have in $Z_k\cap \{u=v_\epsilon=0\}$
$$\chi H(x).\nabla d_\eta= \chi H_2(x)
 \partial_{x_2} d_\eta= {1\over \eta} \chi H_2(x)\chi_{
\{\overline{k}<x_2<\overline{k}+ \eta\}}\geqslant 0.$$

\n \n Let $J=\{w\in(w_1,w_2)\,\,/\,\,\phi_h(w)>t_{\overline{k}}(w)\,\}$.
Then given that $u\in C^{0,1}_{loc}(U)$, one has for some positive constant $C$

\begin{eqnarray*}
| I_3^{\delta\eta} | &\leqslant& \displaystyle{{C\over \eta}
    \int_{Z_k\cap\{u>v_\epsilon=0\}\cap \{\overline{k} < x_2<\overline{k} + \eta\}}
(1- F_\delta (u)) }dx \\
&=& \displaystyle{{C\over \eta} \int_{J}\int_{t_{\overline{k}}(w)}
^{\min (\phi_h(w), t_{\overline{k}+\eta}(w))}(1- F_\delta (u o T_h)) (t,w)
   . (-Y_h(t,w)) }dt dw \\
&\leqslant& C \overline{h}\displaystyle{\int_{J} \Big({1\over \eta} \int_{t_{\overline{k}}(w)}
^{t_{ \overline{x_2}}(w) +{\eta\over \underline{h}}}
  (1- F_\delta (u o T_h)) dt \Big) dw}
\end{eqnarray*}

\n Since the function $t\mapsto 1- F_\delta (u o
T_h(t,w))$ is continuous, we obtain

$$\displaystyle{\limsup_{\eta\rightarrow 0} } | I_3^{\delta\eta}| \leqslant  C \int_{J}
(1- F_\delta (uo T_h(t_{ \overline{k}}(w),w)))
dw.$$ Hence

\begin{eqnarray*}
&&\int_{Z_k^{k+\epsilon}(w_1,w_2)\cap\{0<u-v_\epsilon< \delta\}} {1\over\delta }
\Big( {a(|\nabla u|)\over |\nabla u|} \nabla u-  {a(|\nabla
v_\epsilon|)\over |\nabla v_\epsilon|} \nabla v_\epsilon \Big). \nabla (u-v_\epsilon)^+
 dx\\
&&\quad\leqslant  C \int_{J}(1- F_\delta (uo T_h(t_{\bar k}(w),w))) dw.
\end{eqnarray*}

\n The Lemma follows by letting $\delta\rightarrow 0$.
\qed

\vs 0,5cm

\begin{lemma}\label{l2.4}
Assume that the assumptions of Lemma 2.3 hold. Then we have
\begin{equation}\label{2.10}
\int_{Z_k^{k+\epsilon}(w_1,w_2)}   \mathbb{A}(x)\nabla (u-v_\epsilon)^+ . \nabla\zeta dx\, = \,0 \qquad
\forall\zeta \in \mathcal{D}(Z_k^{k+\epsilon}(w_1,w_2)),
\end{equation}
where $\displaystyle{\mathbb{A}(\xi)=\left(\mathbb{A}_{ij}\right)}$,
$\displaystyle{\mathbb{A}_{ij}={\partial\mathcal{A}^i\over \partial x_j}   }$,
and $\displaystyle{\mathcal{A}^i(\xi)={ a(|\xi|)\over |\xi|} \xi_i}$.
\end{lemma}

\vs0,3cm\n\emph{Proof.} First, we observe that we have for any $\zeta \in \mathcal{D}(Z_k^{k+\epsilon}(w_1,w_2))$
\begin{eqnarray}\label{e2.11}
&& \int_{Z_k^{k+\epsilon}(w_1,w_2)} \chi(\{u>v_\epsilon\}) \Big( {a(|\nabla u|)\over |\nabla
u|}\nabla u -{a(|\nabla v_\epsilon|)\over |\nabla v_\epsilon|}\nabla v_\epsilon \Big). \nabla\zeta dx \nonumber\\
&&=\lim_{\delta\rightarrow 0}
\int_{Z_k^{k+\epsilon}(w_1,w_2)}F_\delta (u-v_\epsilon) \Big( {a(|\nabla u|)\over |\nabla
u|}\nabla u -{a(|\nabla v_\epsilon|)\over |\nabla v_\epsilon|}\nabla v_\epsilon \Big) . \nabla\zeta dx  = \lim_{\delta\rightarrow 0} I_\delta.
\end{eqnarray}

\n Next we have
\begin{eqnarray}\label{e2.12}
 &&I_\delta =  \int_{Z_k^{k+\epsilon}(w_1,w_2)}  \Big( {a(|\nabla u|)\over |\nabla u|}\nabla u
-{a(|\nabla v_\epsilon|)\over |\nabla v_\epsilon|}\nabla v_\epsilon \Big) .  \nabla (F_\delta (u-v_\epsilon)\zeta)dx\nonumber\\
 &&\quad-\,{1\over \delta}
 \int_{Z_k^{k+\epsilon}(w_1,w_2) \cap [0<u-v_\epsilon<\delta]} \zeta  \Big( {a(|\nabla u|)\over |\nabla u|}\nabla u
-{a(|\nabla v_\epsilon|)\over |\nabla v_\epsilon|}\nabla v_\epsilon \Big)  . \nabla  (u-v_\epsilon) dx\nonumber\\
&&\quad= I^1_\delta - I^2_\delta.
\end{eqnarray}

\n By Lemma 2.3, we have
\begin{equation}\label{e2.13}
\lim_{\delta\rightarrow 0}
I_\delta^2=0.
\end{equation}

\n Regarding $I_\delta^1$,  we have from $(P)ii)$ and (2.2), since $(F_\delta (u-v_\epsilon)\zeta
 )\in W^{1,A}_0 (Z_k^{k+\epsilon}(w_1,w_2))$ and $\chi=1$ a.e. in $\{u>v_\epsilon\}$

\begin{eqnarray}\label{e2.14}
&&I^1_\delta= \int_{Z_k^{k+\epsilon}(w_1,w_2)} {a(|\nabla u|)\over |\nabla u|}\nabla u .
\nabla (F_\delta (u-v_\epsilon)\zeta)dx - \int_{Z_k^{k+\epsilon}(w_1,w_2)}
{a(|\nabla v_\epsilon|)\over |\nabla v_\epsilon|}\nabla v_\epsilon. \nabla(F_\delta(u-v_\epsilon)\zeta)dx\nonumber\\
&&\quad =- \int_{Z_k^{k+\epsilon}(w_1,w_2)} \chi H(x).\nabla(F_\delta (u-v_\epsilon)\zeta) dx
+\int_{Z_k^{k+\epsilon}(w_1,w_2)} H(x).\nabla(F_\delta (u-v_\epsilon)\zeta) dx\nonumber\\
&&\quad=0
\end{eqnarray}

\n It follows from (2.11)-(2.14) that
$$\int_{Z_k^{k+\epsilon}(w_1,w_2)} \chi(\{u>v_\epsilon\})  \Big( {a(|\nabla u|)\over |\nabla u|}\nabla u
-{a(|\nabla v_\epsilon|)\over |\nabla v_\epsilon|}\nabla v_\epsilon \Big)  . \nabla\zeta dx\, = \,0 \qquad
\forall\zeta \in \mathcal{D}(Z_k^{k+\epsilon}(w_1,w_2))$$

\n
 which can be written as
\begin{equation}\label{2.15}
\int_{Z_k^{k+\epsilon}(w_1,w_2)}  \chi(\{u>v_\epsilon\})\Big( \int_{0}^{1} { d\over dt} (\mathcal{A} (\nabla w_t)) dt\Big)  . \nabla\zeta dx\, = \,0 \qquad
\forall\zeta \in \mathcal{D}(Z_k^{k+\epsilon}(w_1,w_2))
\end{equation}

\n where $\displaystyle{ \mathcal{A} (\xi)= (\mathcal{A}^1,
\mathcal{A}^2 )(\xi)={ a(|\xi|)\over |\xi|} \xi}$ and $w_t= tu+(1-t)v_\epsilon$.

\n Now observe that

\begin{equation}\label{2.16}
\int_{0}^{1} { d\over dt}
(\mathcal{A} (\nabla w_t)) dt  = \Big(
\int_{0}^{1}{\partial\mathcal{A}^i\over
\partial x_j}(\nabla w_t)\Big)_{i,j=1,2}\nabla (u-v)= \mathbb{A}(x)\nabla (u-v).
\end{equation}

\n Hence we obtain (2.10) from (2.15) and (2.16).
\qed

\vs0,3cm

\begin{lemma}\label{l2.5}

\begin{eqnarray}\label{2.17}
\min(1, a_0) \displaystyle{a(z)\over z}|\xi|^2 \leqslant A_{ij}(z)\xi_i
\xi_j\leqslant  \max(1,a_1) \displaystyle{a(z)\over z}|\xi|^2\qquad\forall z\neq 0,
~~~\forall\xi\in \mathbb{R}^2.
\end{eqnarray}

\end{lemma}

\vs 0,3cm\n\emph{Proof.} Let $z\neq 0$ and $\xi=(\xi_1,\xi_2)\in \mathbb{R}^2$.
Since $A\in C^1(\mathbb{R}^2\setminus\{(0,0)\})$, we have by direct calculation
\begin{eqnarray*}
&&A_{ij}(z)= \displaystyle{\partial
(A(z))_{i}\over\partial z_j}= \Big( \displaystyle{ a'(z) z -
a(z)  \over z^3}\Big) z_i z_j +\displaystyle{a(z)\over z}\delta_{ij}\\
&&A_{ij}(z)\xi_i \xi_j=\Big(
{a'(z) z-a(z)\over z^3}\Big) (z_1\xi_1+z_2\xi_2)^2 + {a(z)\over z} |\xi|^2
\end{eqnarray*}

\n Using (0.1), we obtain

\begin{eqnarray*}
&\displaystyle{a(z)\over z}\Big((a_0-1)\displaystyle{|z.\xi|^2\over z^2}
+ |\xi|^2\Big)\leqslant A_{ij}(z)\xi_i \xi_j\leqslant
\displaystyle{a(z)\over z}\Big((a_1-1)\displaystyle{|z.\xi|^2\over z^2}
+ |\xi|^2\Big) .\end{eqnarray*}

\n Then, if $a_0\geqslant 1$, the left hand side of inequality (2.17) holds.
When $a_0< 1$, we use  Cauchy-Schwartz inequality : $|z.\xi|\leqslant
|z| |\xi|$, to conclude. We proceed in the same way for the right
hand side.\qed

\vs0,3cm

\begin{lemma}\label{l3.4}
Assume that the assumptions of Lemma 2.3 hold. Then only one of
the following situations holds
\begin{eqnarray*}
&&(i) \quad u>0\quad\text{in}\quad Z_k^{k+\epsilon}(w_1,w_2)\\
&&(ii)\quad u=0\quad\text{in}\quad Z_{k+\epsilon}.
\end{eqnarray*}

\end{lemma}

\vs0,3cm\n\emph{Proof.} Assume that $(i)$ is false. Then
\begin{equation*}
\exists(t_0,w_0)\text{ such that }\quad
T_h(t_0,w_0)\in Z_k^{k+\epsilon}(w_1,w_2)\quad\text{ and }\quad u\circ T_h(t_0,w_0)=0.
\end{equation*}
This leads by Theorem 1.1 $ii)$ to
\begin{equation}\label{2.18}
u\circ T_h(t,w_0)=0\quad \forall t\in[t_0,t_{k+\epsilon}].
\end{equation}

\n We will show in this case that $(ii)$ holds. From Lemmas 2.4 and 2.5 we know that

\begin{equation}\label{2.19}
div( \mathbb{A}(x) \nabla (u-v_\epsilon)^+)= 0  \quad\text{in } \quad
Z_k^{k+\epsilon}(w_1,w_2).
\end{equation}

\n Moreover, by Lemma 2.5, the matrix $\mathbb{A}(x)$ satisfies 

\begin{eqnarray}\label{2.20}
&& \min(1, a_0) \lambda(x) |\xi|^2 \leqslant \mathbb{A}(x) \xi.\xi
\leqslant \max(1,a_1) \lambda(x) |\xi|^2\qquad \forall x\in Z_k^{k+\epsilon}(w_1,w_2)~~~\forall \xi\in
\mathbb{R}^2 \nonumber\\
&& \mathrm{with}\qquad  \lambda(x) =\int_0^1 {{a(|\nabla w_t(x)|)}\over{|\nabla w_t(x)|}},
\quad w_t=tu+(1-t)v_\epsilon.
\end{eqnarray}

\n Next, we have $v_\epsilon\in C^{1, \alpha}(Z_k^{k+\epsilon}(w_1,w_2)\cup L)$, where
$$L=\partial Z_k^{k+\epsilon}(w_1,w_2)\cap \{x_2=k+\epsilon\}.$$
We also have $v_\epsilon=0$ on $L$ and $v_\epsilon>0$ in $Z_k^{k+\epsilon}(w_1,w_2)$. So $v_\epsilon$
achieves its minimum value on the line segment $L$. By Lemma 2.2 of \cite{[CL9]}, 
we must have $|\nabla v_\epsilon|>0$ along $L$.
Therefore for $\delta$ small enough such that $w_1+\delta<w_2-\delta$,
there exists two positive constants $c_0$, $c_1$ such that
\begin{equation}\label{e2.22}
c_0 \leqslant |\nabla v_\epsilon(x)| \leqslant c_1\quad
\forall x\in \overline{Z}_k^{k+\epsilon}(w_1,w_2)\cap\{k+\epsilon-\delta\leqslant x_2\leqslant k+\epsilon\}
\cap \{w_1+\delta\leqslant w\leqslant w_2-\delta\}=Z_{k+\epsilon-\delta}^{k+\epsilon}.
\end{equation}

\n On the other hand,  $|\nabla  u|$ is also bounded in
$Z_{k+\epsilon-\delta}^{k+\epsilon}$ since
$u\in C^{0,1}(\overline{Z}_k^{k+\epsilon}(w_1,w_2))$ (\cite{[CL9]} Theorem 1.1). It follows from (2.20)-(2.21) that we
have for two positive constants $\lambda_0$ and $\lambda_1$
\begin{equation*}
\lambda_0\leqslant\lambda(x)\leqslant \lambda_1 \qquad \mathrm{in
}\quad Z_{k+\epsilon-\delta}^{k+\epsilon}
\end{equation*}

\n and therefore we get from (2.20)
\begin{equation}\label{e2.23}
\min(1, a_0) \lambda_0 |\xi|^2 \leqslant \mathbb{A}(x) \xi.\xi
\leqslant \max(1,a_1) \lambda_1 |\xi|^2\qquad 
\forall x\in Z_{k+\epsilon-\delta}^{k+\epsilon}~~~\forall \xi\in
\mathbb{R}^2
\end{equation}

\n Taking into account (2.18), we see that
\begin{equation}\label{2.22}
Z_{k+\epsilon-\delta}^{k+\epsilon}\cap \{u=0\}\neq \emptyset.
\end{equation}

\n It follows from (2.19), (2.22), (2.23), and the strong maximum principle
that $(u-v_\epsilon)^+\equiv 0$ in $Z_{k+\epsilon-\delta}^{k+\epsilon}$.
Consequently, we obtain  $u\leqslant v_\epsilon$ in 
$Z_{k+\epsilon-\delta}^{k+\epsilon}$, 
and then $u oT_h(t_{k+\epsilon}(w),w ) =0$ for all $w\in ( w_1+\delta,w_2-\delta)$.
Since $\delta$ is arbitrary small, we get $u oT_h(t_{k+\epsilon}(w),w ) =0$
for all $w\in ( w_1,w_2)$. Hence $(ii)$ holds by Theorem 1.1 $ii)$.

\qed
\vs 0,3cm
\begin{lemma}\label{l2.7} Let $w_0 \in(w_*, w^*)$, $x_0=T_h(t_0,w_0)$
such that for some $\eta>0$, $B_\eta(T_h(t_0,w_0))\subset \subset U$
and $u(x_0)=0$. Then there exists two sequences $(t_n^-,w_n^-)_n$
and $(t_n^+,w_n^+)_n$ such that for all $n$:
\begin{eqnarray*}
& i)~ T_h(t_n^-,w_n^-)\in B_\eta(T_h(t_0,w_0))\cap \{w<w_0\},\quad u\circ T_h(t_n^-,w_n^-)=0,
\quad \displaystyle{\lim_{n\rightarrow\infty}(t_n^-,w_n^-)=(t_0,w_0)}.\\
& ii)~ T_h(t_n^+,w_n^+)\in B_\eta(T_h(t_0,w_0))\cap \{w>w_0\},\quad u\circ T_h(t_n^+,w_n^+)=0,
\quad \displaystyle{\lim_{n\rightarrow\infty}(t_n^+,w_n^+)=(t_0,w_0)}.
\end{eqnarray*}

\end{lemma}

\vs 0,3cm\n\emph{Proof.} First we observe that by Lemma 1.6 the
following situations cannot hold simultaneously
\begin{eqnarray*}
& a)~ u\circ T_h>0\mbox{ in } B_\eta(T_h(t_0,w_0))\cap \{w<w_0\}\\
& b)~ u\circ T_h>0\mbox{ in } B_\eta(T_h(t_0,w_0))\cap \{w>w_0\}.
\end{eqnarray*}

\n In fact, to prove the lemma, it is enough to show that neither $a)$ nor $b)$ holds.
So assume for example that $a)$ holds. Then by Lemma 1.6 there exists a sequence
$(t_n^+,w_n^+)_n$ such that $T_h(t_n^+,w_n^+)\in B_\eta(T_h(t_0,w_0))\cap \{w>w_0\}$,
\[u\circ T_h(t_n^+,w_n^+)=0\quad
\mbox{ and } \quad\displaystyle{\lim_{n\rightarrow\infty}(t_n^+,w_n^+)=(t_0,w_0)}.\]

\vs 0.2cm \n Let $k=\max\{T_h^2(t_0,w_0),T_h^2(t_n^+,w_n^+)\}$. 
Then since $u(x_0)=0$ and $u$  is continuous at $x_0$, we may assume that for
$n$ large enough, we have
\begin{equation}\label{e2.24}
   u\circ T_h(t_k(w),w) \leqslant \vartheta_{\epsilon}(\epsilon) \qquad
   \forall w\in (w_0,w_n^+).
\end{equation}

\vs 0,2cm\n For $\epsilon>0$ small enough and $n$ large enough, we may assume that

\begin{equation}\label{e2.25}
Z_k^{k+\epsilon}(w_0,w_n^+)\subset\subset U
\end{equation}

\vs 0,2cm\n Using (2.24), (2.25) and Lemma 2.6, we conclude that 
for $\epsilon>0$ small enough and $n$ large enough, we have $u=0$ in 
$Z_{k+\epsilon}\cap  T_h(\{w_0< w<w_n^+\})$.
Now since we have assumed that $a)$ holds, we are in contradiction with Lemma 1.6.

\vs 0,2cm\n Similarly, if we assume that $b)$ holds, we will get a contradiction
as well.

\qed

\vs 0,5cm\n Finally, by using Lemma 2.6 and Lemma 2.7, we can
establish the main result of the paper.

\vs 0,3cm
\begin{theorem}\label{l2.1} The
function $\phi_h$ is continuous in the interval $(w_*,w^*)$.
\end{theorem}

\vs 0.3cm \n \emph{Proof}. Let $w_0 \in(w_*,w^*)$. We will prove that
$\phi_h$ is continuous at $w_0$. To this end, it is enough to show that 
$\phi_h$ is upper semi-continuous at $w_0$.

\vs 0.2cm \n Let $x_0=T_h(\phi_h(w_0), w_0)=T_h(t_0, w_0)$ and let $\epsilon>0$.
Since $u(x_0)=0$ and $u$  is continuous at $x_0$,  there exists
$\eta \in (0,\epsilon)$ such that
  \begin{equation}\label{e2.26}
   u(x) \leqslant \vartheta_{\epsilon}(\epsilon) \qquad
   \forall x\in B_{\eta} (x_0)\subset\subset U.
\end{equation}

\vs 0.2cm \n By Lemma 2.7 there exists two sequences $(t_n^-,w_n^-)_n$
and $(t_n^+,w_n^+)_n$ such that for all $n$:
\begin{eqnarray*}
& i)~ T_h(t_n^-,w_n^-)\in B_\eta(T_h(t_0,w_0))\cap \{w<w_0\},\quad u\circ T_h(t_n^-,w_n^-)=0,
\quad \displaystyle{\lim_{n\rightarrow\infty}(t_n^-,w_n^-)=(t_0,w_0)}.\\
& ii)~ T_h(t_n^+,w_n^+)\in B_\eta(T_h(t_0,w_0))\cap \{w>w_0\},\quad u\circ T_h(t_n^+,w_n^+)=0,
\quad \displaystyle{\lim_{n\rightarrow\infty}(t_n^+,w_n^+)=(t_0,w_0)}.
\end{eqnarray*}

\vs 0.2cm \n Let $k=\max\{T_h^2(t_n^-,w_n^-),T_h^2(t_0,w_0),T_h^2(t_n^+,w_n^+)\}$
and let $C$ be the constant in Lemma 1.2.
We observe that we can choose $\epsilon$ small enough and $n$ large enough so that
\begin{eqnarray}\label{e2.27}
&& \epsilon'=\epsilon/2C<\underline{h}/{2\overline{h}}\nonumber\\
&& Z_k^{k+\epsilon'}(w_n^-,w_n^+)\subset\subset U\nonumber\\
&&\partial Z_k^{k+\epsilon'}(w_n^-,w_n^+)\cap\{x_2=k\}\subset B_\eta(x_0).
\end{eqnarray}

\vs 0.2cm \n Using (2.26), (2.27) and Lemma 2.6, we see that for $n$ large enough, we have
$$u=0\qquad \hbox{ in }T_h(\{w_n^-< w<w_n^+\})\cap \{x_2\geqslant k
+\epsilon'\}.$$

\n Therefore we obtain

\begin{equation}\label{2.28}
\phi_h(w) \leqslant t_{k+\epsilon'}( w)\quad \forall w\in (w_n^-,w_n^+).
\end{equation}

\n From Lemma 1.2, we infer that if $\eta<\epsilon/4C$:

\begin{eqnarray}\label{2.29}
t_{k+\epsilon'}( w)&\leqslant& t_{x_{02}}(w_0)+C(|k+\epsilon'-x_{02}|+|w-w_0|)\nonumber\\
&\leqslant& t_0+C(\eta+\epsilon'+\eta)=t_0+2C\eta +\epsilon/2\nonumber\\
&\leqslant& t_0+\epsilon/2+\epsilon/2=t_0+\epsilon.
\end{eqnarray}
\n Combining (2.28) and (2.29), we obtain

$$\phi_h(w) \leqslant \phi_h ( w_0) + \epsilon \quad \forall w\in (w_n^-,w_n^+)$$

\n which is the upper semi-continuity of $\phi_h$ at $w_0$.\qed

\vs 0,5cm \n \begin{remark}
Thanks to Theorem 2.1 and (1.3), and 
taking into account that $\chi=1$ a.e. in $\{u>0\}$, we obtain

\[\chi=\chi_{\{u>0\}}\]
\end{remark}

\vs 0.5cm \n\emph{Acknowledgments } The authors are
grateful for the facilities and excellent research conditions provided
by Fields Institute where part of this research was carried out.

\end{document}